\documentclass[bibother]{asl}
\usepackage{rotate}

\usepackage[UKenglish]{babel}
\usepackage[utf8]{inputenc}
\usepackage{verbatim}

\usepackage{theoremref}
\usepackage{hyperref}

\usepackage{dsfont}
\usepackage{enumitem}

\usepackage{thmtools}

\usepackage{tikz-cd}

\usepackage{catoptions}
\makeatletter

\def\Autoref#1{%
	\begingroup
	\edef\reserved@a{\cpttrimspaces{#1}}%
	\ifcsndefTF{r@#1}{%
		\xaftercsname{\expandafter\testreftype\@fourthoffive}
		{r@\reserved@a}.\\{#1}%
	}{%
		\ref{#1}%
	}%
	\endgroup
}
\def\testreftype#1.#2\\#3{%
	\ifcsndefTF{#1autorefname}{%
		\def\reserved@a##1##2\@nil{%
			\uppercase{\def\ref@name{##1}}%
			\csn@edef{#1autorefname}{\ref@name##2}%
			\autoref{#3}%
		}%
		\reserved@a#1\@nil
	}{%
		\autoref{#3}%
	}%
}

\makeatother

\theoremstyle{plain}
\newtheorem{theorem}{Theorem}[section]
\newtheorem{thmx}{Theorem}

\newtheorem{lemma}[theorem]{Lem\-ma}
\newtheorem{proposition}[theorem]{Prop\-o\-si\-tion}
\newtheorem{question}[theorem]{Question}

\newtheorem{fact}[theorem]{Fact}

\theoremstyle{definition}
\newtheorem{definition}[theorem]{Definition}
\newtheorem{construction}[theorem]{Construction}
\newtheorem{remark}[theorem]{Remark}
\newtheorem{example}[theorem]{Ex\-am\-ple}
\newtheorem{claim}{Claim}

\newcommand{\N}[0]{\mathbb{N}}
\newcommand{\Z}[0]{\mathbb{Z}}
\newcommand{\Q}[0]{\mathbb{Q}}
\newcommand{\R}[0]{\mathbb{R}}

\renewcommand{\H}[0]{\mathbf{H}}
\newcommand{\OO}[0]{\mathcal{O}}
\newcommand{\MM}[0]{\mathcal{M}}

\newcommand{\supp}[0]{\mathrm{supp}}
\newcommand{\brackets}[1]{\left( #1 \right)}
\newcommand{\setbr}[1]{\left\{ #1 \right\}}
\newcommand{\pow}[1]{\!\left(\!\left( #1 \right)\!\right)}

\renewcommand{\L}[0]{\mathcal{L}}
\newcommand{\Lor}[0]{\mathcal{L}_{\mathrm{or}}}
\newcommand{\Log}[0]{\mathcal{L}_{\mathrm{og}}}
\newcommand{\Lr}[0]{\mathcal{L}_{\mathrm{r}}}
\newcommand{\Lvf}[0]{\mathcal{L}_{\mathrm{vf}}}
\newcommand{\Lovf}[0]{\mathcal{L}_{\mathrm{ovf}}}

\renewcommand{\div}[1]{{#1}^{\mathrm{div}}}

\newcommand\restr[2]{{
		\left.\kern-\nulldelimiterspace 
		#1
		\vphantom{\big|} 
		\right|_{#2}
}}
\newcommand{\vmin}[0]{v_{\min}}

\newcommand{\one}[0]{\mathds{1}}

\newcommand{\h}[1]{{#1}^{\mathrm{h}}}

\DeclareMathOperator{\Th}{Th}

\usepackage{xcolor}
\usepackage[normalem]{ulem}
\usepackage{cancel}

\def\confname{}
\def\copyrightyear{2000}

\title[Definability of henselian valuations]{Definability of henselian valuations by conditions on the value group}

\author[L.~S.~Krapp]{Lothar Sebastian Krapp}
\revauthor{Krapp, Lothar Sebastian}
\author[S.~Kuhlmann]{Salma Kuhlmann}
\revauthor{Kuhlmann, Salma}
\author[M.~Link]{Moritz Link}
\revauthor{Link, Moritz}

\address{Fachbereich Mathematik und Statistik\\Universität Konstanz\\78457 Konstanz, Germany}
\email{sebastian.krapp@uni-konstanz.de}
\urladdr{http://www.math.uni-konstanz.de/\texttildelow krapp/}
\thanks{The first author was partially supported by Werner und Erika Messmer Stiftung.}

\address{Fachbereich Mathematik und Statistik\\Universität Konstanz\\78457 Konstanz, Germany}
\email{salma.kuhlmann@uni-konstanz.de}
\urladdr{https://www.mathematik.uni-konstanz.de/kuhlmann/}

\address{Fachbereich Mathematik und Statistik\\Universität Konstanz\\78457 Konstanz, Germany}
\email{moritz.link@uni-konstanz.de}
\urladdr{https://www.mathematik.uni-konstanz.de/volkwein/members/research-group/\\doctoral-students/moritzlink/}
\thanks{During parts of this project, the second and third authors were supported by the German Excellence Strategy and the German Research Foundation.}

\begin{document}	
	
	\begin{abstract}
		Given a henselian valuation, we study its definability (with and without parameters) by examining conditions on the value group.
		We show that any henselian valuation whose value group is not closed in its divisible hull is definable in the language of rings, using one parameter. Thereby we strengthen known definability results.
		Moreover, we show that in this case, one parameter
		is optimal in the sense that one cannot obtain definability without parameters.
		To this end, we present a construction method for a $t$-henselian non-henselian ordered field elementarily equivalent to a henselian field with a specified value group.
	\end{abstract}
	
	\maketitle
	
	\section{Introduction}
	
	In recent years, the study of definable henselian valuations has received a considerable amount of attention and several definability conditions have been established. Early results of this study were related to decidability questions (cf.~e.g.~\cite{robinson2,ax}), whereas a more recent motivation is due to the Shelah--Hasson Conjecture on the classification of NIP fields (cf.~e.g.~\cite{shelah,jahnke,johnson1,dupont,halevi3,halevi2,halevi,krapp2}). We refer the reader to \cite{fehm} for a more detailed survey on the definability of henselian valuations. While some of these results deal with the existence of non-trivial definable henselian valuations and others consider the quantifier complexity of defining formulas (cf.~\cite[Sections~2~\&~3]{fehm}), in this paper we are interested in the definability of a given henselian valuation with a certain value group.\footnote{There are also conditions for the definability of a given henselian valuation only depending on its residue field (see our \nameref{sec:conlude}).} 
	
	The strongest currently known definability results which only pose conditions on the value group of a given henselian valuation are exhibited in \cite{hong} and \cite{krapp}. In the following, we state the results under consideration for this paper and outline how we strengthen these. 
	
	\begin{fact}{\rm\cite[Corollary~2]{hong}.}\thlabel{fact:hongdiscrete}
		Let $(K,v)$ be a henselian valued field such that the value group $vK$ is discretely ordered. Then $v$ is $\Lr$-definable with one parameter from $K$.
	\end{fact}
	
	For densely ordered value groups, the following is established.
	
	\begin{fact}{\rm\cite[Theorem~3]{hong}.}\thlabel{fact:hongprendible}
		Let $(K,v)$ be a henselian valued field such that $vK$ is densely ordered. Assume that $vK$ contains a convex subgroup that is $p$-regular but not $p$-divisible, for some prime $p\in \N$. Then $v$ is $\Lr$-definable with one parameter from $K$.
	\end{fact}
	
	While in \cite{hong} definability is treated in the language of rings $\Lr$, the definability results in \cite{krapp} are stated for ordered fields in the language of ordered rings $\Lor$.
	
	\begin{fact}{\rm\cite[Theorem~5.3~(2)]{krapp}.}\thlabel{fact:kkl}
		Let $(K,<,v)$ be an ordered henselian valued field such that $vK$ is not closed in its divisible hull. Then $v$ is $\Lor$-definable  with one parameter from $K$.
	\end{fact}
	
	After introducing preliminary notions in \Autoref{sec:prelim}, we show 	
	in \Autoref{sec:rendible} that the hypothesis on $vK$ in \Autoref{fact:hongprendible} is strictly stronger than that in \Autoref{fact:kkl} (see \Autoref{prop:prendibleclosed} and \Autoref{ex:counterhong}). Thus, 	 { we obtain the following strengthening of both \Autoref{fact:hongprendible} and \Autoref{fact:kkl}:\footnote{Theorem~\ref{thm:A} will be restated as \Autoref{thm:strenhong}.}}

	\begin{thmx} \label{thm:A}
		Let $(K,v)$ be a henselian valued field such that $vK$ is not closed in its divisible hull. Then $v$ is $\Lr$-definable with one parameter from $K$.
	\end{thmx}

	In \Autoref{sec:0def}, we examine definability \emph{without parameters} of the valuations. 
	 {To this end, in \Autoref{constr:main} -- the technical heart of this paper -- we refine \cite[page~338]{prestelziegler} to obtain a $t$-henselian non-henselian field (see \Autoref{def:thens})
	with prescribed value groups. More precisely, given an arbitrary family of non-trivial ordered abelian groups $\{A_{\gamma}\mid \gamma \in -\omega\}$, \Autoref{constr:main} produces a $t$-henselian non-henselian ordered field $(K,<)$ as well as, for each $n\in\N$,  a convex valuation $u_n$ on $K$  whose value group is the Hahn sum (see page~\pageref{hahnsum}) over $\{A_{\gamma}\mid \gamma \leq -n\}$.}
	
	 {By applying \Autoref{constr:main} to the specific case where all $A_{\gamma}$ are identical, we obtain the following:\footnote{Theorem~\ref{thm:B} will be restated as \Autoref{thm:secondmain}.}}
	
	\begin{thmx}\label{thm:B}
		Let $A\neq\{0\}$ be an ordered abelian group. Then there exists an ordered henselian valued field $(L,<,v)$ such that $vL$ is the Hahn sum $\coprod_{-\omega}A$ and the only parameter-free $\Lor$-definable henselian valuation on $(L,<)$ is the trivial valuation. In particular, $v$ is not $\Lor$-definable without parameters.
	\end{thmx}
	
	 {Finally, by choosing appropriately the group in Theorem~\ref{thm:B}, we can show that} \Autoref{fact:hongdiscrete} -- \ref{fact:kkl} and Theorem~\ref{thm:A} are optimal in the sense that,  {in general}, one cannot obtain definability of the valuations \emph{without} parameters (see \Autoref{ex:counter1} and \Autoref{ex:counter2}).
	
	
	
	We conclude this work in \Autoref{sec:questions} by collecting open questions motivated by our results.
	
	 {}
	
	\section{Preliminaries}\label{sec:prelim}
	
	More details on the model and valuation theoretic notions we use can be found in \cite{marker,engler,kuhlmann}.
	
	We denote by $\N$ the set of natural numbers \emph{without} $0$ and by $\omega$ the set of natural numbers \emph{with} $0$. All topological properties are considered with respect to the order topology. Whenever an algebraic structure carries a standard linear ordering,  we refer to this specific ordering if not explicitly stated otherwise.
	
	Let $K$ be a field and let $v$ be a valuation on $K$. We denote the \textbf{valuation ring} of $v$ (i.e.\ the subring $\{a\in K\mid v(a)\geq 0\}$ of $K$) by $\OO_v$ and its \textbf{valuation ideal} (i.e.\ its maximal ideal $\{a\in \OO_v\mid v(a)>0\}$) by $\MM_v$. Moreover, we denote the \textbf{value group} of $v$ by $vK$ and the \textbf{residue field} $\OO_v/\MM_v$ by $Kv$. An element $a+\MM_v \in Kv$ with $a\in \OO_v$ is also denoted by $av$. Suppose that $K$ admits an ordering $<$. Then we say that $v$ is \textbf{convex} (with respect to $<$) if $\OO_v$ is a convex subset of $K$ (with respect to the ordering $<$). By \cite[Lemma~2.1]{knebusch}, any henselian valuation on an ordered field is convex.
	Given a valued field $K$ we denote by $\h{K}$ its \textbf{henselisation}. 	
	
	Let $\Log=\{+,-,0,<\}$ be the \textbf{language of ordered groups}, $\Lr=\{+,-,\cdot,0,\allowbreak 1\}$ the \textbf{language of rings} and  $\Lor=\Lr\cup\{<\}$ the \textbf{language of ordered rings}. If no confusion is likely to arise, we denote an $\Log$-structure $(G,+,-,0,<)$ simply by $G$, an $\Lr$-structure $(R,+,-,\cdot,0,1)$  by $R$ and an $\Lor$-structure $(R,+,-,\cdot,0,1,\allowbreak<)$ by $(R,<)$.  
	We say that a valuation $v$ on a field $K$ (an ordered field $(K,<)$) is \textbf{$\Lr$-definable} (\textbf{$\Lor$-definable}) if $\OO_v$ is an $\Lr$-definable  ($\Lor$-definable) subset of $K$. Here, definability is meant \emph{with} parameters. We say that $v$ is \textbf{$0$-$\Lr$-definable} (\textbf{$0$-$\Lor$-definable}) if there is a defining formula without parameters from $K$.
	
	If we expand $\Lr$ by a unary predicate $\OO_v$, we obtain the \textbf{language of valued fields} $\Lvf$. Similarly, $\Lovf=\Lor\cup\{\OO_v\}$ denotes the \textbf{language of ordered valued fields}. Let $(K,\OO_v)$ be a valued field. An atomic formula of the form $v(t_1)\geq v(t_2)$, where $t_1$ and $t_2$ are $\Lr$-terms, is equivalent to $t_1=t_2=0\vee (t_2\neq 0\wedge \tfrac{t_1}{t_2}\in \OO_v)$. Hence, by abuse of notation, we can also write $(K,v)$ for the $\Lvf$-structure of a valued field. Likewise the $\Lovf$-structure of an ordered valued field is denoted by $(K,<,v)$. 
	
	
	Let $(\Gamma,<)$ be a linearly ordered set, and for each $\gamma\in \Gamma$, let $A_\gamma\neq\{0\}$ be an abelian group (expressed additively). The corresponding \textbf{Hahn product} over $\left\{A_\gamma\mid \gamma\in \Gamma\right\}$ is given by
	$$\underset{\gamma \in \Gamma}{\H} A_\gamma=\setbr{\left.\!s\in \prod_{\gamma \in \Gamma}A_\gamma\ \right|\  \supp(s)\text{ is well-ordered}}\!,$$
	where $\prod_{\gamma \in \Gamma}A_\gamma$ denotes the group product over $(A_\gamma)_{\gamma\in \Gamma}$ and $\supp(s)$ denotes the \textbf{support} of $s$ given by the set $\{\gamma\in \Gamma\mid s(\gamma)\neq 0\}$. The Hahn product is an  abelian group under pointwise addition. It carries a valuation with value set $\Gamma$ given by $\vmin\colon s\mapsto \min\supp(s)$	for $s\neq 0$. 
	Any $s\in \H_{\gamma\in \Gamma}A_\gamma$ can be expressed as a sum $s=\sum_{\gamma \in \Gamma}s_\gamma \one_\gamma$, where $s_\gamma = s(\gamma)$ and $\one_\gamma$ is the characteristic function mapping $\gamma$ to some fixed $1_\gamma \in A_\gamma$ (which in our application is always clear from the context) and everything else to $0$.
	The valued subgroup $$\coprod_{\gamma \in \Gamma}A_\gamma=\{s\in \H_{\gamma \in \Gamma}A_\gamma \mid \supp(s) \text{ is finite}\}$$
	is called the \textbf{Hahn sum} \label{hahnsum} over $(\Gamma,(A_\gamma)_{\gamma\in \Gamma})$.
	If $\Gamma$ is a convex subset of $\Z$ with $\min \Gamma=a$ and $\max \Gamma =b$ (where we set $a=-\infty$ if $\Gamma$ is not bounded from below, and $b=\infty$ if $\Gamma$ is not bounded from above), then we express $\coprod_{\gamma \in \Gamma} A_\gamma$ as $\coprod_{\gamma=a}^b A_\gamma$ if $a\neq - \infty$, and as $\coprod_{\gamma > - \infty}^b A_\gamma$ if $a = -\infty$.
	For an ordinal $\alpha$ we have that $\coprod_{\gamma\in-\alpha}A_\gamma = \H_{\gamma\in-\alpha}A_\gamma$.
	If there is an ordered abelian group $A$ such that for each $\gamma \in \Gamma$ we have $A_\gamma=A$, we simply write $\H_\Gamma A$ ($\coprod_{\Gamma} A$) for the Hahn product (Hahn sum) above.
	 {If $A_\gamma$ is an ordered abelian group for each $\gamma \in \Gamma$, then also $\H_{\gamma\in\Gamma}A_\gamma$ is an ordered abelian group with the ordering given as follows: for any non-zero $s=\sum_{\gamma\in\Gamma} s_\gamma \one_\gamma \in \H_{\gamma\in\Gamma} A_\gamma$ we have  $s > 0$ if and only if $s_{\gamma_0} > 0$, where $\gamma_0 = \vmin \left(s\right)$.}

	Let $k$ be a field and let $G$ be an ordered abelian group. The \textbf{Hahn field} (also called \textbf{power series field}) $k\pow{G}$ has as underlying additive group the Hahn product $\H_{G}k$. Multiplication on $k\pow{G}$ is defined by
	$(rs)_h=\sum_{g\in G}r_g s_{h-g}$ for any $h\in G$. 
	To distinguish between elements of Hahn products and Hahn fields, we write an element $s\in k\pow{G}$ as sum $s=\sum_{g\in G}s_gt^g$, where $t^g$ is the characteristic function with $t^g(g)=1$ and $t^g(h)=0$ for any $h\neq g$. 
	The valuation $\vmin$ on $k\pow{G}$ is henselian.
	Again, if $k$ is ordered, then also $k\pow{G}$ carries a natural ordering as described above making it an ordered field. We denote by $k[G]$ the \textbf{group ring} with coefficient field $k$ and value group $G$. It consists of all elements of $k\pow{G}$ with finite support. Its quotient field, denoted by $k(G)$, is generated in $k\pow{G}$ by the set of all monic monomials, i.e.\ $k(G)=k(t^g\mid g\in G)$. 
	
	Let $G$ be an ordered abelian group. We say that $G$ is \textbf{discretely ordered} if there exists some $1\in G$ with $0<1$ and there is no element in $G$ strictly between $0$ and $1$. If $G$ is not discretely ordered, then it is \textbf{densely ordered}. Let $n\in \N$. Then $G$ is called \textbf{$n$-regular} if every infinite convex subset of $G$ contains at least one $n$-divisible element. 	If $G$ is $n$-regular for every $n\in \N$, then we say that it is \textbf{regular}.
	We denote the \textbf{divisible hull} of $G$ by $\div{G}$. 
	
	
	\section{Definable valuations}\label{sec:rendible}
	
	We start by stating and proving the main result of this section. Subsequently, we will discuss how this generalises the known definability results from \Autoref{fact:hongprendible} and \Autoref{fact:kkl}. 
	
	\begin{theorem} \thlabel{thm:strenhong}
		Let $(K,v)$ be a henselian valued field. Suppose that $vK$ is not closed in $\div{vK}$\!. Then $v$ is $\Lr$-definable with one parameter.
	\end{theorem}
	
	\Autoref{thm:strenhong} can essentially be proven by the methods in \cite[page~15]{hong} and \cite[page~19]{hong2}.\footnote{We thank Blaise Boissonneau for pointing out this method to the first author.} However, due to the more general context there are certain adjustments to be made, which is why we present the proof in detail.
	
	\begin{proof}	
		We set $G=vK$. Let $\varepsilon \in K$ and let $n\in \N$ such that $v(\varepsilon)$ is not $n$-divisible in $G$ but $\tfrac{v(\varepsilon)}{n}$ is a limit point of $G$ in $\div{G}$. Since $\tfrac{v(\varepsilon)}{n}$ is a limit point of $G$,  also   $-\tfrac{v(\varepsilon)}{n}$ is a limit point of $G$. Thus, by replacing $\varepsilon$ by $\varepsilon^{-1}$ if necessary, we may assume that $\tfrac{v(\varepsilon)}{n}$ is a left-sided limit point of $G$, i.e.\ for any positive $h\in \div{G}$ the half-open interval $\left (\tfrac{v(\varepsilon)}{n}-h,\tfrac{v(\varepsilon)}{n}\right ]$ in $\div{G}$ contains some element from $G$.
		
		Set $\Phi_\varepsilon=\{x\in K\mid v(\varepsilon x^n)>0\}$. We  show that the formula $\varphi(x,\varepsilon)$ given by
		$$\exists y\ (y^n-y^{n-1}=\varepsilon x^n)$$
		defines $\Phi_\varepsilon$.		
		Let $x\in \Phi_\varepsilon$. Then $\varepsilon x^n\in \MM_v$ and by Hensel's Lemma there exists $y\in K$ such that $y^n-y^{n-1}=\varepsilon x^n$, i.e.\ $K \models \varphi(x,\varepsilon)$.
		Conversely, let $x\in K$ such that $K \models \varphi(x,\varepsilon)$, i.e.\ there exists $y\in K$ such that $y^{n-1}(y-1)=\varepsilon x^n$. By applying $v$, we obtain
		$$(n-1)v(y)+v(y-1) = v(\varepsilon)+nv(x)$$
		and distinguish three cases. 
		
		\emph{Case 1:} $v(y)>0$. Then $v(y-1)=0$ and thus 
		$0<(n-1)v(y)=v(\varepsilon)+nv(x)=v(\varepsilon x^n)$, as required.
		
		\emph{Case 2:} $v(y)=0$. Then $v(y-1)=v(\varepsilon)+nv(x)$. If $v(y-1)>0$, then we are done immediately. Otherwise, $v(y-1)=0$ and we obtain  $\tfrac{v(\varepsilon)}{n} = - v(x)$, a contradiction to the fact that $v(\varepsilon)$ is not $n$-divisible in $G$.
		
		\emph{Case 3:} $v(y)<0$. Then we obtain $v(\varepsilon) = n (v(y)-v(x))$, again a contradiction, as $v(\varepsilon)$ is not $n$-divisible.
		
		Now set
		$$\Psi_\varepsilon = \{\varepsilon x^n\mid x\in \Phi_\varepsilon\}=\{\varepsilon x^n\mid x\in K, v(\varepsilon x^n)>0\}.$$
		Again, $\Psi_\varepsilon$ is $\Lr$-definable only using the parameter $\varepsilon$.
		Finally, set 
		$$\Omega_\varepsilon=\setbr{x^n-x^{n-1}\mid x\in K, K\models\exists y\ \exists (z\in \Psi_\varepsilon)\  z(y^n-y^{n-1})=x^n-x^{n-1} }\!.$$
		We show that $\Omega_\varepsilon =\MM_v$, as then $\OO_v$ is defined by $\forall (u\in \Omega_\varepsilon) \ xu\in \Omega_\varepsilon$.
		
		First let $a\in\Omega_\varepsilon$ and fix $x\in K$ such that $x^n-x^{n-1} = a$. Then by definition of $\Omega_\varepsilon$ there are $y\in K$ and $z\in \Psi_\varepsilon$ such that $z(y^n-y^{n-1})=a$. Let $r\in K$ with $z=\varepsilon r^n\in \MM_v$. Assume, for a contradiction, that $v(a)\leq 0$. If $v(a)<0$, then 
		\begin{align*}
		    0 &>  v(x^{n-1}(x-1))=(n-1)v(x)+v(x-1)=nv(x)\\
		    \text{and} \quad 0 &> v(z)+(n-1)v(y)+v(y-1)>(n-1)v(y)+v(y-1)=nv(y).
		\end{align*}
		Since $$\varepsilon=zr^{-n}=(x^n-x^{n-1})[(y^n-y^{n-1})r^n]^{-1},$$ we obtain
		$$v(\varepsilon)=v(x^n-x^{n-1})-v(y^n-y^{n-1})-nv(r)=n(v(x)-v(y)-v(r)),$$
		contradicting the choice of $\varepsilon$. If $v(a)=0$, then 
		$0=(n-1)v(x)+v(x-1)$ and thus $v(x)=v(x-1)=0$. Moreover, $0=v(z)+(n-1)v(y)+v(y-1)$ and thus $0>-v(z)=(n-1)v(y)+v(y-1)=nv(y)$. This implies 
		$$v(\varepsilon)=v(x^n-x^{n-1})-v(y^n-y^{n-1})-nv(r)=n(v(y)-v(r)),$$
		which again contradicts the fact that $v(\varepsilon)$ is not $n$-divisible.
		This shows that $v(a)>0$, i.e.\ $a\in \MM_v$.
		
		For the converse, let $a\in \MM_v$. By Hensel's Lemma there exists $x\in K$ such that $x^n-x^{n-1}=a$. Since  $\tfrac{v(\varepsilon)}{n}$ is a left-sided limit point of $G$, we obtain some $b \in K^\times$ with
		$$\frac{v(\varepsilon)}{n}-\frac{v(a)}{n}<v(b)< \frac{v(\varepsilon)}{n}.$$
		Hence, $-v(\varepsilon)<-nv(b)<v(a)-v(\varepsilon)$. We set $z=\varepsilon b^{-n}$ and obtain that $0<v(z)<v(a)$, whence $z\in \Psi_\varepsilon$ and $v\!\brackets{\frac{a}{z}}>0$. Again by applying Hensel's Lemma, we find some $y\in K$ such that $y^n-y^{n-1}=\frac{a}{z}$. This yields $a\in \Omega_\varepsilon$ and therefore $\Omega_\varepsilon = \MM_v$, as required.
	\end{proof}
	
	
	
	\begin{remark}
	    Note that the class of non-trivial ordered abelian groups that are closed in their divisible hull is elementary in the language of ordered groups. Indeed, it is axiomatised by the axiom scheme
	    \begin{equation*}
        \forall a \:\left( \forall \left(b>0\right)\: \exists g \:\: |a-ng|<b \rightarrow \exists h \:\: a=nh \right),
        \end{equation*}
        one for each $n\in\N$.\qed
	\end{remark}
	
	Our next aim is to show that \Autoref{thm:strenhong} implies \Autoref{fact:hongprendible}. Note that in both of these statements the value group is densely ordered {; in \Autoref{fact:hongprendible} by assumption and in \Autoref{thm:strenhong} since any discretely ordered abelian group is already closed in its divisible hull.}
	
	\begin{proposition}
		\thlabel{prop:prendibleclosed}\thlabel{cor:prendibleclosed}
		Let $G$ be a densely ordered abelian group. Suppose that there exists some  $n\in \N$ such that $G$ contains a convex subgroup $C$ that is $n$-regular but not $n$-divisible. Then $G$ is not closed in $\div{G}$.	
	\end{proposition}
	
	\begin{proof}
		Note that since $G$ contains a convex subgroup that is not $n$-divisible, we have that $G$ is a non-divisible ordered abelian group. Assume, for a contradiction, that $G$ is closed in $\div{G}$. Since $C$ is not $n$-divisible, there exists $a \in C$ such that $\tfrac{a}{n} \in \div{G}\setminus G$. Moreover, since $G$ is closed in $\div{G}$, there exists a positive $\delta \in \div{G}$ such that
		$$\left(\frac{a}{n}-\delta,\frac{a}{n}+\delta\right) \cap G =\emptyset.$$
		As $0$ is a limit point of $G$ in $\div{G}$ (cf.\ \cite[Lemma~3.6]{krapp}), there exists $b\in G$ with $0<b<\delta$ and hence
		\begin{equation}
		\label{eq:closed}
		\left(\frac{a}{n}-\frac{b}{n},\frac{a}{n}+\frac{b}{n}\right) \cap G =\emptyset.
		\end{equation}
		By convexity of $C$, we may choose $b\in C^{>0}$. Further, since $G$ is densely ordered and thus $C$ is densely ordered, the interval $(a-b,a+b)$ is an infinite convex subset of $C$. By $n$-regularity of $C$, there is some $z\in C$ such that $a-b \leq nz \leq a+b$.
		This is equivalent to
		$$\frac{a-b}{n} \leq z \leq \frac{a+b}{n},$$
		a contradiction to (\ref{eq:closed}). Hence, $G$ is not closed in $\div{G}$.
	\end{proof}
	
	
	
	\Autoref{cor:prendibleclosed} shows that the condition on the value group in \Autoref{fact:hongprendible} is covered by the condition in \Autoref{thm:strenhong}.
	The converse of \Autoref{prop:prendibleclosed} does not hold, as verified by the following example.
	
	\begin{example}
		\thlabel{ex:counterhong}
		We consider the ordered abelian group 
		\begin{equation*}
			G = \left(\coprod_{\omega} \Z\right) + a\Z\subseteq \underset{\omega}{\H}\ \Z, \text{ where }a=\sum_{n\in \omega} 2\one_n.
		\end{equation*}
		 {Note that $G$ consists of all elements $s = \sum_{n\in\omega} s_n \one_n \in \H_{\omega} \Z$ for which the following holds: there exist $m\in\Z$ and $n_0\in\omega$ such that for any $n \geq n_0$ we have $s_n = 2m$. We now show that $G$ is not closed in $\div{G}$ but for any $n\in\N$ and any convex subgroup $C$ of $G$, if $C$ is $n$-regular, then it is already $n$-divisible.} 
		
		By \cite[Lemma~3.6]{krapp}, we have that $G$ is densely ordered, as its value set $\vmin G=\omega$ does not have a maximum.
		We first note that the divisible hull of $G$ is given by 
		\begin{equation*}
			\div{G}= \left(\coprod_{\omega} \Q\right) + a\Q.
		\end{equation*}
		Let $s = \sum_{n\in \omega} \one_n \in \div{G}$.  {By the description of the elements of $G$ above, one immediately sees that $s \notin G$. We now show that $s$ is a limit point of $G$ in $\div{G}$, whence $G$ is not closed in $\div{G}$}.  {Let $b\in \div{G}$ be positive and set $\ell = \vmin (b)$. Further, let $g = \sum_{n=0}^{\ell } \one_n \in G$ and note that $g<s$.  As $b_\ell > 0$, we obtain that $0<s - g = 0 \cdot \one_\ell + \sum_{n\geq \ell+1}\one_n < b$. We have thus found an element $g\in G$ with $s-b<g<s$, showing that $s$ is a limit point of $G$ in $\div{G}$.}
		
		Let $n\in \N$ with $n\geq 2$.	We show that $G$ contains no non-trivial $n$-regular convex subgroup. Thus, in particular, the only $n$-regular convex subgroup of $G$ is the trivial group, which is $n$-divisible. 
		
		Let $C$ be a non-trivial convex subgroup of $G$.
		Then there exists $\ell \in \omega$ such that $C$ is of the form
		\begin{equation*}
			C = \{g\in G\mid \vmin(g)\geq \ell\}= \left(\coprod_{\gamma=\ell}^\infty \Z\right) + \setbr{\left.\! \sum_{i\geq \ell} k\one_i\ \right|\ k\in 2\Z}\!
		\end{equation*}
		(cf.~\cite[page~50~f.]{kuhlmann}).
		Now let $c,d \in C$ with $c={q}\one_{\ell}$ and $d={q}\one_{\ell}+\one_{\ell+1}$, 
		where ${q}\in \N$ is some prime with ${q}>n$. Then $c<d$ and $(c,d)$ is an infinite convex subset of the densely ordered group $C$. Since for any $z \in (c,d)$ we have that $z_\ell = q$ and $n \nmid q$ in $\Z$, there is no $n$-divisible element in $(c,d)$. Hence, $C$ is not $n$-regular.\qed
	\end{example}
	
	\Autoref{ex:counterhong} shows that the condition in \Autoref{thm:strenhong} is indeed more general than the condition in \Autoref{fact:hongprendible}. More precisely, let $k$ be any field and let $G$ be the ordered abelian group constructed in 	\Autoref{ex:counterhong}; then \Autoref{thm:strenhong} implies that the henselian valuation $\vmin$ on $k\pow{G}$ is $\Lr$-definable (with one parameter), whereas  \Autoref{fact:hongprendible} cannot be applied directly to obtain the $\Lr$-definability of $\vmin$.
	
	In conclusion, \Autoref{thm:strenhong} is a strict generalisation of both  \Autoref{fact:hongprendible} and \Autoref{fact:kkl} for densely ordered value groups: \Autoref{thm:strenhong} clearly generalises \linebreak \Autoref{fact:kkl}, as it has the same condition on the value group but does not rely on the existence of an ordering on the field. Moreover, 
	\Autoref{cor:prendibleclosed} shows that the value group condition of \Autoref{thm:strenhong} implies that of \Autoref{fact:hongprendible} but \Autoref{ex:counterhong} illustrates that the converse does not hold.
	
	\section{ {Main construction and }parameter-free definability}\label{sec:0def}
	The main aim of this section is to show that the conclusions of \Autoref{fact:hongdiscrete} and \Autoref{thm:strenhong}  {can, in general,} not be strength\-ened to $0$-$\Lr$-definability. More precisely, we construct a class of ordered henselian valued fields that do not admit any non-trivial $0$-$\Lor$-definable henselian valuation (see \Autoref{thm:secondmain}). This class contains ordered henselian valued fields $(K,<,v)$ with  discretely ordered value group $vK$ (see \Autoref{ex:counter1}) as well as such with a value group $vK$ not closed in $\div{vK}$ (see \Autoref{ex:counter2}).

	The method for constructing an \emph{ordered} henselian valued field $(K,<,v)$ that does not admit a non-trivial $0$-$\Lor$-definable henselian valuation relies on the construction of a $t$-henselian non-henselian \emph{ordered} field.\footnote{Our construction is inspired by the arguments in \cite[page~338]{prestelziegler}, which have also been refined in \cite[Section~6]{fehm2}, \cite[Section~6]{jahnke5} and \cite[Section~4]{anscombe}.} 
	We start by clarifying the notion of $t$-henselianity, which was first introduced in \cite[Section~7]{prestelziegler}.\footnote{Several characterisations for a field to be $t$-henselian are given in \cite[page~354]{fehm2}.}
	
	\begin{definition}\thlabel{def:thens}\thlabel{lem:thens}
	    We call a field \textbf{henselian} if it admits a non-trivial henselian valuation.
		A field $K$ is \textbf{$t$-henselian} if there exists a henselian field $L$ such that $K\equiv L$ (in the language $\Lr$).
	\end{definition}
	
	We now establish a sufficient condition for $t$-henselianity based on $n_\leq$-henselian valuations: 
	Let $n\in \N$. Following the terminology of \cite[Definition~6.1]{fehm2}, we say that a valued field $(K,v)$ is \textbf{$n_\leq$-henselian} if Hensel's Lemma holds in $(K,v)$ for all polynomials of degree at most $n$. Note that $(K,v)$ is henselian if and only if it is $n_\leq$-henselian for any $n\in \N$. Arguing similarly as in the proof of \cite[Theorem~4.1.3]{engler}, one can show that $(K,v)$ is $n_\leq$-henselian if and only if every polynomial $X^\ell+X^{\ell-1}+a_{\ell-2}X^{\ell-2}+\ldots+a_0 \in \OO_v[X]$ with $\ell\leq n$ and $a_{\ell-2},\ldots,a_0\in \MM_v$ has a zero in $K$. Moreover, by the arguments in the proof of  \cite[Lemma~2.1]{knebusch}, for $n\geq 2$ any $n_\leq$-henselian valuation on a real field $K$ is convex with respect to any ordering on $K$.

	
	

	
	\begin{lemma}
	\thlabel{lem:nthenselian}
		Let $K$ be a field. Suppose that for any $n\in \N$, there exists a non-trivial $n_\leq$-henselian valuation on $K$. Then $K$ is $t$-henselian.
	\end{lemma}
	
	\begin{proof}
		Let $h_n$ be the $\Lvf$-sentence stating that $v$ is a non-trivial valuation for which Hensel's Lemma holds for polynomials of degree at most $n$. Then any finite subset of the $\Lvf$-theory $\Th(K)\cup\{h_n\mid n\in \N\}$ is satisfied by $(K,v_m)$ for some $m\in \N$, where $v_m$ is $m_\leq$-henselian. Hence, by the Compactness Theorem, there exists a henselian field $L$ with $L\models \Th(K)$, i.e.\ $L\equiv K$.
	\end{proof}
	

	
	Our aim is to refine the construction of a $t$-henselian non-henselian field in \cite[page~338]{prestelziegler} by gaining influence on the resulting value group. More precisely, starting with a family of non-trivial ordered abelian groups $\{A_{\gamma}\mid \gamma\in -\omega\}$, we construct a non-henselian ordered field $(K,<)$ such that for any $n\in \N$, there is an $n_\leq$-henselian valuation $v_n$ on $K$ with	$v_nK = \coprod_{\gamma>-\infty}^{-n} A_{\gamma}$.
	
	
	In the following, for an ordered abelian group $G$ we consider $\div{G}$ also as a $\Q$-vector space with the usual scalar multiplication. For an algebraic field extension $K\subseteq L$ and an element $\alpha \in L$, we denote by $\deg_K(\alpha)$ the degree of the minimal polynomial of $\alpha$ over $K$. 
	
    \Autoref{cor:biljakovic} and \Autoref{lem:nonhens} below will be used in \Autoref{constr:main}. Note that the statement of \Autoref{lem:nonhens} can be adjusted for an arbitrary field $k$ of characteristic $0$ (instead of an ordered field).
	
	\begin{lemma}
		\thlabel{cor:biljakovic}
		Let $k$ be a field of characteristic $0$ and let $G$ be an ordered abelian group. Moreover, let $a\in k[G]$ such that for some $n\geq 2$ we have $\supp(a) = \{g_1,\dots,g_n,0\}$ and $g_1,\dots,g_n$ are linearly independent in the $\Q$-vector space $\div{G}$. Then $a$ is prime in $k[G]$.
	\end{lemma}
	
	\begin{proof}
		First note that the group of units of $k[G]$ is given by the multiplicative group of monomials of $k\pow{G}$, i.e.\
		$k[G]^\times=\{ct^g\mid c\in k^\times, g\in G\}$. 
		By \cite[Proposition~4.1]{biljakovic}, every irreducible element of $k[G]$ is prime. It thus suffices to show that $a$ is irreducible. Let $r,s\in k[G]$ with $a=rs$. We obtain from \cite[Proposition~4.8]{biljakovic} that $a$ is prime in $k[\div{G}]$ and hence irreducible in $k[\div{G}]$ (as $k[\div{G}]$ is an integral domain). We may thus assume that $r$ is a unit in $k[\div{G}]$. Hence, $r\in k[G]\cap \{ct^g\mid c\in k^\times, g\in \div{G}\}=k[G]^\times$. This shows that $r$ is a unit in $k[G]$, as required.
	\end{proof}
	
	For two ordered abelian groups $A$ and $B$, we denote by $A\amalg B$ the abelian group $A\times B$ ordered lexicographically.
	
	\begin{lemma} 
		\thlabel{lem:nonhens}
		Let $(k,<)$ be an ordered field and let $A_1, A_2\neq\{0\}$ be ordered abelian groups. 
		Moreover, set $K=k(A_1\amalg A_2)$ and let $(\h{K},<)$ be a henselisation of $K$ in $(k\pow{A_1\amalg A_2},<,\vmin)$. Then for any $n\in \N$, there exists an intermediate field extension $K\subseteq L \subsetneq \h{K}$ such that $(L,\vmin)$ is $n_\leq$-henselian but not henselian.
	\end{lemma}
	
	\begin{proof}
		Let $n\in \N$ and let $p\in \N$ be prime with $p>n$.
		Since $p_{\leq}$-henselianity implies $n_\leq$-henselianity, it suffices to construct an ordered valued field $(L,\linebreak <,\vmin)$ with $(K,<,\vmin)\subseteq (L,  <,\vmin) \subsetneq (\h{K},<,\vmin)$
		such that for any polynomial of the form
		\begin{align*}
			f(T) = T^m + a_{m-1}T^{m-1}+ \ldots + a_0 \in L[T],
		\end{align*}
		where $m\leq p$, and any $\alpha\in \h{K}$ with $f(\alpha)=0$, we have that $\alpha\in L$.
		Let $M$ be the set of all natural numbers $m$ such that for any prime number $q\in \N$ dividing $m$ we have $q\leq p$. 
		We set $L=K\!\left(\alpha\in \h{K}\mid \deg_K(\alpha)\in M\right)$, i.e.\ $L$ is generated by all elements in $\h{K}$ whose degree over $K$ only has prime factors less than or equal to $p$. Note that since $M$ is closed under multiplication and divisors, $\deg_K(x)\in M$ for any $x\in L$.
		Let $m\leq p$ and let $f(T)=T^m + a_{m-1}T^{m-1}+ \ldots + a_0  \in L[T]$. By construction of $L$, we have $\deg_K(a_i)\in M$ for any $i\in\{0,\ldots,m-1\}$. Moreover, since $[K(a_0,\ldots,a_{m-1}):K]$ divides the product $[K(a_0):K]\ldots[K(a_{m-1}):K]=\deg_K(a_0)\ldots \deg_K(a_{m-1})$, we obtain $[K(a_0,\ldots,a_{m-1}):K] \in M$. Suppose that $\alpha\in \h{K}$ is a zero of $f$. Since also $[K(a_0,\ldots,a_{m-1},\alpha):K(a_0,\ldots,a_{m-1})]\leq m\leq p$, we obtain  
		\begin{align*}&[K(a_0,\ldots,a_{m-1},\alpha):K]=\\ &\qquad[K(a_0,\ldots,a_{m-1},\alpha):K(a_0,\ldots,a_{m-1})]\cdot  [K(a_0,\ldots,a_{m-1}):K]\in M.\end{align*}
		Finally, since $\deg_K(\alpha)=[K(\alpha):K]$ divides $[K(a_0,\ldots,a_{m-1},\alpha):K]$, we obtain $\deg_K(\alpha)\in M$ and thus $\alpha \in L$.
		
		It remains to show that $(L,\vmin)$ is not henselian, i.e.\ that there exists $y\in \h{K}$ such that $y\notin L$. Let $q \in \N$ be a prime number greater than $p$. We set
		$$ a = t^{(b_1,0)} + t^{(0,b_2)} + 1 \in k[A_1\amalg A_2],$$
		{ where $b_1\in A_1$ and $b_2\in A_2$ are fixed positive elements,
		and observe that $a$ is prime in $k[A_1\amalg A_2]$ by \Autoref{cor:biljakovic}. Further, let $$g(T) = T^q + a T^{q-1} + a t^{(b_1,0)} T^{q-2} + \dots + a t^{(b_1,0)} T + a t^{(b_1,0)} \in K[T].$$
		Note that $\vmin(a)=0$ and $\vmin(a t^{(b_1,0)})=(b_1,0)>0$.
		Hence, by applying the residue map, we obtain the polynomial $g\vmin (T)=T^q+T^{q-1}\in k[T]$. Since $g\vmin(T)$ has a simple zero in $k$, we obtain that $g$ has a zero $y\in \h{K}$. 
		Moreover, by Eisenstein's Criterion we obtain that $g$ is irreducible over $K$ and thus 
		$q = [K(y) : K]$.
	}Since $q$ does not divide any element in $M$, we obtain by construction of $L$ that $y\notin L$. 
		Hence, $L\subsetneq \h{K}$ as required.
	\end{proof}
	
	We now proceed with our main construction.
	
	\begin{construction}[Main Construction]\thlabel{constr:main}
		Let  {$\{A_\gamma\mid \gamma\in \omega\}$ be a family of non-trivial ordered abelian groups and set $G:= \coprod_{\gamma\in-\omega}A_{-\gamma}$.} We construct an ordered field $(K,<)$ such that for any $n\in \N$ there exists a convex $n_\leq$-henselian valuation $ u_n$ on $K$ with value group $u_nK = \coprod_{\gamma > -\infty}^{-n} A_ {-\gamma}$ and $K$ does not admit a non-trivial henselian valuation. In particular, $K$ is $t$-henselian (by \Autoref{lem:nthenselian}) but not henselian. 
		

		By induction on $m\in\omega$ we construct a chain of ordered fields $(K_m, <)$ each endowed with a valuation $v_\mathrm{min}^m$. From this the ordered field $(K,<)$ will be obtained by an inverse limit construction, as well as a chain of convex valuations $v_m$ on $K$ for $m\in\omega$. Ultimately, for each $n\in \omega$ we set $u_{2n}$ to be $v_{n}$ and $u_{2n+1}$ to be a certain coarsening of $u_{2n}$.
		
		 {To initiate the induction, we} set $(K_{0},<)=(\R,<)$. Now let $m\in \N$ and suppose that the ordered field $(K_{m-1},<)$ has already been defined.
		We denote the valuation given by $s\mapsto \min \supp (s)$ on  {$K_{m-1}\pow{A_{2m-1}\amalg A_{2(m-1)}}$} by $\vmin^m$.
		By \Autoref{lem:nonhens}, there is an ordered valued field $K_m$ with
		\begin{align*} 
		    &K_{m-1}\left(A_{2m-1}\amalg A_{2(m-1)}\right) \subseteq K_m \\
		    \subsetneq\: &\h{K_{m-1}\left(A_{2m-1}\amalg A_{2(m-1)}\right)}\subseteq K_{m-1}\pow{A_{2m-1}\amalg A_{2(m-1)}},
		\end{align*}
		where the ordering is inherited from  {$K_{m-1}\pow{A_{2m-1}\amalg A_{2(m-1)}}$} and the valuation is $\vmin^m$, 
		such that any zero in  {$\h{K_{m-1}\!\left(A_{2m-1}\amalg A_{2(m-1)}\right)}$} of a polynomial over $K_m$ of degree at most $ 2m+1$ already lies in $K_m$.
		Since the extensions above are immediate, we have $K_m\vmin^m =K_{m-1}$. Let $\OO_m'$ and $\MM_m'$ be the valuation ring and valuation ideal of $\vmin^m$ in $K_m$ and let $\pi_m\colon \OO_m' \to K_{m-1}$ denote the corresponding residue map $a \mapsto a\vmin^m$.
		
		For $i,j\in \omega$ with $j\leq i$ we define the map $ \psi_{i,j}\colon K_i \cup \{\infty\} \to K_j \cup \{\infty\}$ as follows:
		We set
		$$\psi_{i,i} = \mathrm{id}_{K_i \cup \{\infty\}} \text{ and }\psi_{i,i-1} (x) = \begin{cases}
		\pi_i (x)& \text{ if } x\in \OO_{i}',\\
		\infty& \text{ if } x\notin\OO_{i}'.
		\end{cases}$$
		For $0\leq j<i-1$ we set 
		$\psi_{i,j} = \psi_{j+1,j} \circ \dots \circ \psi_{i,i-1}$.
		Now we define $I$ to be the inverse limit of $\left(K_m\cup \{\infty\}\right)_m$, i.e.\ $I$ is the following set of sequences:
		$$I = \left\{\left.\! (x_m) \in \prod_{m\in\omega} \left(K_m\cup \{\infty\}\right)\ \right\vert\ \psi_{i,j}(x_i)=x_j \text{ for all } i,j\in \omega \text{ with }j\leq i \right\}\!.$$
		 {Note that if for a given sequence $(x_m) \in I$, we have $x_\ell=0$ (respectively $x_\ell=\infty$) for some $\ell\in\omega$, then $x_k=0$ (respectively $x_k=\infty$) for any $k\leq \ell$. On the other hand, if $x_\ell \notin\{0,\infty\}$ for some $\ell\in \omega$, then 
		\begin{align}\label{eq:idealeq}
		   x_k \in \OO^\prime_k \setminus \MM^\prime_k,
		\end{align}
		 i.e.\ $\vmin^k(x_k)=0$, for all $k>\ell$.
		Therefore,  if $x_\ell \notin\{0,\infty\}$ and $x_{\ell-1} \in\{0,\infty\}$  for some $\ell\in \N$, then $\vmin^\ell(x_\ell) \neq 0.
		$}
		\begin{claim}\label{cl:Kfield}
		The set $K:= I\setminus\{(\infty)\}$ can be made an ordered field.
		\end{claim}
		\begin{proof}[Proof of Claim~\ref{cl:Kfield}]
		\renewcommand{\qedsymbol}{$\Diamond$}
		Addition, multiplication and an ordering on $K$ can be defined by a standard procedure\footnote{See \cite[page~114\,ff.]{engler} for analogous construction in the group case.} making $K$ an ordered field: Let $(x_m),(y_m) \in K$ and let $\ell \in \omega$ be least such that $x_\ell\neq \infty$ and $y_\ell\neq \infty$. For any $i\geq \ell$ set $z_i=x_i+y_i$. For any $i<\ell$ set $z_i=\psi_{\ell,i}(x_\ell+y_\ell)$. Then $(x_m) + (y_m)$ is given by $(z_m) \in K$. Likewise, multiplication can be defined via the multiplication operations on the fields $K_m$. Finally, a non-zero element $(x_m)$ of $K$ is positive if and only if $x_i$ is non-negative for any $i\in \omega$ with $x_i\neq \infty$.
		\end{proof}
		
		For any $n\in\omega$ we define the projection map
		$\psi_n \colon K \to K_n\cup\{\infty\}, 
		(x_m) \mapsto x_n$. This gives rise to the following commutative diagram and equation for $j\leq i$:
		\begin{center}
		\begin{tikzcd}
           K \arrow[rd, "\psi_j"'] \arrow[r, "\psi_i"] & K_i \cup \{\infty\} \arrow[d, "{\psi_{i,j}}"] \\
            & K_j \cup \{\infty\}            
        \end{tikzcd}
        \begin{align}\label{eq:places}
            \psi_{i,j} \circ \psi_{i} = \psi_j
        \end{align}
        \end{center}
		One can easily verify that $\psi_n$ is a place, and therefore $\OO_n = \psi_n^{-1}(K_n)$ defines a valuation ring on $K$ with maximal ideal $\MM_n=\psi_n^{-1}(\{0\})$ and residue field $K_n$. By \eqref{eq:places} we compute for any $a\in \MM_{n+1}$: $$\psi_n(a) = \psi_{n+1,n}(\psi_{n+1}(a)) = \psi_{n+1,n}(0) = 0.$$
		Therefore we have $a\in \MM_n$, and thus $\OO_n \subseteq \OO_{n+1}$.
		
		We denote the valuation on $K$ with valuation ring $\OO_n$ by $v_n$ and fix $n\in\N$.  
		\begin{claim}\label{cl:vnlesshens}
		$v_n$ is a $( 2n+1)_\leq$-henselian valuation on $K$.
		\end{claim}
		\noindent We observe that, if \Autoref{cl:vnlesshens} is established, then we also have that $v_n$ is $( 2n)_\leq$-henselian and convex (as $ 2n+1\geq 2$).
		\begin{proof}[Proof of Claim~\ref{cl:vnlesshens}]
		\renewcommand{\qedsymbol}{$\Diamond$}
		Let 
		$$f(T) = T^{ 2n+1} + T^{ 2n} + a^{( 2n-1)}T^{ 2n-1} + \dots + a^{(0)} \in K[T],$$
		where $a^{( 2n-1)},\dots,a^{(0)}\in \MM_n$. We need to find a zero $(y_m)\in K$ of $f$.
		For any $i\in\omega$, set $a_i^{( 2n-1)} = \psi_i\!\left(a^{( 2n-1)}\right),\dots,a_i^{(0)}=\psi_i\!\left(a^{(0)}\right) \in K_i$ and 
		$$f_i(T) = T^{ 2n+1} + T^{ 2n} + a_i^{( 2n-1)}T^{ 2n-1} + \dots + a_i^{(0)}\in K_i[T].$$
		First let $i\leq n$. Then $a_i^{( 2n-2)}= \dots = a_i^{(0)}=0$
		and therefore $f_i(T) = T^{ 2n+1} + T^{ 2n}$. Set $y_i=-1$, which is a simple zero of $f_i$. 
		Now let $i>n$. We construct $y_i$ iteratively. To this end, assume that $y_{i-1}\in K_{i-1}$ is a simple zero of $f_{i-1}$ with $\psi_{i-1,j}(y_{i-1})=y_j$ for $0\leq j\leq i-1$.
		Applying $\pi_i$ to the coefficients of $f_i(T)\in K_i[T]$ we obtain the residue polynomial
		$$T^{ 2n+1} + T^{ 2n} + a_{i-1}^{( 2n-1)} T^{ 2n-1} + \dots + a_{i-1}^{(0)}=f_{i-1}(T)\in K_{i-1}[T].$$
		Thus, $f_i$ has a simple zero $y_i$ in the henselisation of $K_i$ with $\pi_i(y_i)=y_{i-1}$. By construction of $K_i$ (as $i>n$), we obtain $y_i\in K_i$, as required.
		Finally, one can verify by direct computation that the element $(y_m)\in K$ we obtained from this construction is a zero of $f$.
		
		The argument above shows that $v_n$ satisfies the condition for $( 2n+1)_\leq$-hen\-se\-li\-an\-i\-ty only for polynomials of degree equal to $ 2n+1$. However, since the coarsening of a $( 2n+1)_\leq$-henselian valuation is also $( 2n+1)_\leq$-henselian (see the proof of \cite[Lemma~A.4.24]{krappdiss}), we immediately obtain $( 2n+1)_\leq$-henselianity of $v_n$.
		\end{proof}
		
		It remains to show that $v_n$ is not henselian.
		\begin{claim}
		\label{cl:vnnothens}
		$v_n$ is not henselian.
		\end{claim}
		\begin{proof}[Proof of Claim~\ref{cl:vnnothens}]
	    \renewcommand{\qedsymbol}{$\Diamond$}
		Assume for a contradiction that $v_n$ is henselian. Let $m\in\N$ and let
		$$
		g(X)=X^m + X^{m-1} + a_{m-2}X^{m-2} + \dots + a_0 \in K_{n+1}[X]
		$$
		for some $a_{m-2},\dots,a_0\in \MM'_{n+1}$. Let $a^{(m-2)},\dots,a^{(0)}\in K$ such that $\psi_{n+1}\!\left(a^{(i)}\right)=a_i$ for $0\leq i \leq m-2$. Then $a^{(m-2)},\dots,a^{(0)}\in \psi_n^{-1}(\{0\})=\MM_n$. As $v_n$ is henselian, there exists $x\in \OO_n$ such that
		$x^m + x^{m-1} + a^{(m-2)}x^{m-2} + \dots + a^{(0)} = 0$.
		A direct computation shows that $\psi_{n+1}(x)$ is a zero of the polynomial
		$g(X)$.
		Thus, we have that $\left(K_{n+1},\vmin^{n+1}\right)$ is henselian, which is a contradiction to the construction of $K_{n+1}$.
		\end{proof}
		
		So far we have seen that for any $n\in\N$ there exists a non-henselian but $( 2n+1)_\leq$-henselian valuation $v_n$ on $K$. It remains to show that $K$ does not admit a non-trivial henselian valuation.
		
		\begin{claim}
		\label{cl:nohens}
		$K$ does not admit a non-trivial henselian valuation.
		\end{claim}
		\begin{proof}[Proof of Claim~\ref{cl:nohens}]
		\renewcommand{\qedsymbol}{$\Diamond$}
		Let $u$ be any non-trivial convex valuation on $K$ and pick $a=(a_m)\in K$ with $u(a)<0$, i.e.\ $a\notin \OO_u$. Let $\ell\in \N$ be least such that $a_\ell\neq \infty$. Then $\psi_\ell(a)=a_{\ell}\in  K_\ell$, whence $a\in \OO_\ell$. Since $v_\ell$ is convex, we obtain that $v_\ell$ is a proper coarsening of $u$. Therefore, since $v_\ell$ is not henselian, also $u$ is not henselian. Finally, as any henselian valuation on an ordered field is convex (cf.~\cite[Lemma~2.1]{knebusch}), this suffices to show that $K$ does not admit a non-trivial henselian valuation.
		\end{proof}
		
		Finally, we show inductively that $v_nK = \coprod_{\gamma>-\infty}^{-2n}A_{-\gamma}$. Since $v_n$ is a coarsening of $v_0$ for all $n \in \N$, we have that $v_nK$ is the quotient of $v_0K$ by $v_0(\OO^\times_n)$. Therefore our next task is to compute $v_0K$ and $v_0(\OO^\times_n)$.  {In order to do so, let us first recall the construction of the fields $K_m$.}
		
		We have the inclusions as ordered fields
		$$ \R(A_1\amalg A_0)\subseteq K_1\subseteq \R\pow{A_1\amalg A_0}\subseteq \R\pow{G},$$
		where the last inclusion is obtained by embedding  {$A_1 \amalg A_0$ into $G$ as its convex subgroup $\coprod_{\gamma = -1}^0 A_{-\gamma}$.}
		In the next step of our iterative field construction we obtained the ordered field inclusions
		 {$$\R(A_1\amalg A_0)(A_3\amalg A_2)\subseteq K_2\subseteq \R\pow{A_1\amalg A_0}\pow{A_3\amalg A_2}\subseteq \R\pow{G}.$$}
		The last inclusion is given via the isomorphism\footnote{Generally, for any ordered abelian groups $C$ and $D$ and any ordered field $(k,<)$ we have $(k\pow{C}\pow{D},<)\cong (k\pow{D\amalg C},<)$ via the isomorphism of ordered fields induced by $at^ct^d\mapsto at^{(d,c)}$ for any $a\in k$ and $(c,d)\in C\times D$.}  {$$\R\pow{A_1\amalg A_0}\pow{A_3\amalg A_2}\cong \R\pow{(A_3\amalg A_2)\amalg(A_1\amalg A_0)}$$ and the embedding of $(A_3\amalg A_2)\amalg(A_1\amalg A_0)$ into $G$ as its convex subgroup $\coprod_{\gamma=-3}^0 A_{-\gamma}$.}
		Continuing this process, we have that  {$$\R\left(\coprod_{\gamma=-(2m-1)}^0 A_{-\gamma}\right) \subseteq K_m \subseteq \R\pow{\coprod_{\gamma=-(2m-1)}^0 A_{-\gamma}}\subseteq \R\pow{G}.$$}
		We denote by $\vmin$ the standard valuation $s\mapsto \min \supp (s)$ on $\R\pow{G}$. By the inclusion described above, $\vmin$ restricted to $K_m$ defines a valuation with value group  {$\vmin K_m=\coprod_{\gamma=-(2m-1)}^0 A_{-\gamma}$. For any $m\in\omega$ we have that $\vmin K_m$  can be regarded as a convex subgroup of $G$. Note that the valuation $\vmin$ restricted to $K_m$ is a refinement of $\vmin^m$: Both $\vmin$ and $\vmin^m$ are convex, $ \vmin K_m=\coprod_{\gamma=-(2m-1)}^0 A_{-\gamma}$ and $\vmin^m K_m = \left(A_{2m-1} \amalg A_{2(m-1)}\right)$. Thus, $\vmin^m$ is the coarsening of $\vmin$ induced by the convex subgroup $\coprod_{\gamma=-(2m-3)}^0 A_{-\gamma}$ of $\vmin K_m$. 
		}
		
		 {We now want to define a map $v$ on $K$ such that $v$ is a valuation on $K$ and $\OO_v = \OO_0$. We set $v(0) = \infty$.}
		
        \begin{claim} 
                    \label{cl:valcalc}
                    For any element $(x_m) \in K^\times$ and for any $x_i, x_j \notin \{0,\infty\}$ we have that $\vmin(x_i) = \vmin(x_j)$.
        \end{claim}	
        \begin{proof}[Proof of Claim~\ref{cl:valcalc}]
        \renewcommand{\qedsymbol}{$\Diamond$} 
        Let $\ell\in \omega$ be maximal with $x_{\ell-1} \in \{0,\infty\}$, if such $\ell$ exists, and $\ell=0$ otherwise. It suffices to show that $\vmin(x_k) = \vmin(x_\ell)$ for any $k\geq \ell$. We proceed by induction on $k$. For $k=\ell$ there is nothing to prove. Now assume that $\vmin(x_\ell) = \vmin(x_k)$ for some fixed $k \geq \ell$. We first express $x_k \in K_k$ as
        $$x_k = \sum_{g \in G} c_g t^g \in \R\pow{G}.$$
        Now set $g_0 = \vmin(x_k)$. Then $\vmin(x_k)= \vmin(x_\ell) \in \coprod_{\gamma = -(2\ell-1)}^0 A_{-\gamma} \subseteq G$ if $\ell\geq 1$, and $g_0=\vmin(x_0)=0$ if $\ell=0$. 
        Further we set $$g_0 = \sum_{\gamma=-(2\ell-1)}^0 h_{\gamma} \one_{\gamma},$$
        where $h_{\gamma}\in A_{-\gamma}$ for each $\gamma$.
        By \eqref{eq:idealeq} we have that $\vmin^{k+1}(x_{k+1}) = 0$ and by construction $x_{k+1}\vmin^{k+1} = x_k$. Therefore we can set 
        $$x_{k+1} = c_{g_0} t^{g_0} + \sum_{g>g_0} c_gt^g + \varepsilon,$$
        for some $\varepsilon \in \MM^\prime_{k+1}$. Considering $\varepsilon$ as an element of $\R\pow{G}$ we can write
        $$\varepsilon = \sum_{g\in \coprod^0_{\gamma = -(2k+1)} A_{-\gamma}} \varepsilon_g t^g.$$
        Let $\delta = \vmin(\varepsilon)$ and write
        $$\delta = \delta_{-(2k+1)} \one_{-(2k+1)} + \delta_{-2k} \one_{-2k} + \sum_{\gamma = -(2k-1)}^{-2\ell} \delta_{\gamma} \one_{\gamma} + \sum_{\gamma=-(2\ell-1)}^0 \delta_{\gamma} \one_{\gamma}.$$
        Now since $\varepsilon\in \MM^\prime_{k+1}$, we have that $\delta_{-(2k+1)} \one_{-(2k+1)} + \delta_{-2k} \one_{-2k} > 0$. Hence,
        \begin{align*}
        	\delta - g_0 &= \underbrace{\delta_{-(2k+1)} \one_{-(2k+1)} + \delta_{-2k} \one_{-2k}}_{>0} +\! \sum_{\gamma = -(2k-1)}^{-2\ell} \!\delta_{\gamma} \one_{\gamma} +\! \sum_{\gamma=-(2\ell-1)}^0 \!\left(\delta_{\gamma} - h_{\gamma}\right) \one_{\gamma} \\&> 0,
        \end{align*} 
    	and therefore $\delta > g_0$.
        Thus, $\vmin(x_{k+1}) = \vmin(x_k + \varepsilon) = \vmin(x_k) \allowbreak= \vmin(x_\ell)$, as required.
        \end{proof}
		
		Let $(x_m)\in K^\times$. We define $$v((x_m)):=\vmin(x_i) \in G$$ for any $x_i \notin \{0, \infty\}$.   
		
		\begin{claim}
		\label{cl:vval}
		The map $v$ defines a valuation on $K$ with $vK = G$.
		\end{claim}
		\begin{proof}[Proof of Claim~\ref{cl:vval}]
		\renewcommand{\qedsymbol}{$\Diamond$}
		Note that by Claim~\ref{cl:valcalc}, $v$ is well-defined. Let $(x_m),(y_m)\in K^\times$ and let $i$ be large enough such that $x_i,y_i\notin\{0,\infty\}$. We obtain
		$$v((x_m)(y_m))=\vmin(x_iy_i)=\vmin(x_i)+\vmin(y_i)=v((x_m))+v((y_m)).$$
		Similarly, we prove that $v((x_m) + (y_m)) \geq \min\{v((x_m)), v((y_m))\}.$
		
		Moreover, $vK = G$: For any  $g \in G$ we can choose the least $\ell$ such that there is an element $a_\ell\in K_\ell$ with $\vmin(a_\ell)=g$. Now for any $(x_m) \in K$ with $x_\ell=a_\ell$ a direct calculation shows $v((x_m))=\vmin(a_\ell)=g$.
		\end{proof}
		%
		
		We now show that $\OO_v=\OO_0$, this will yield that $v_0K$ is (isomorphic to) $G$, as required.
		\begin{claim}
		\label{cl:valrings}
		$\OO_v = \OO_0$.
		\end{claim}
		\begin{proof}[Proof of Claim~\ref{cl:valrings}]
		\renewcommand{\qedsymbol}{$\Diamond$}
		Let $(x_m) \in K$. First suppose that $(x_m) \notin \OO_0=\linebreak\psi_0^{-1}(K_0)$. Then $x_0=\infty$. Let $i$ be maximal such that $x_i=\infty$. Then $v((x_m))=\vmin(x_{i+1})<0$. Hence, $(x_m) \notin \OO_v$. For the converse, suppose that $(x_m) \in \OO_0\setminus\{0\}$. Then $x_0\neq \infty$. If $x_0=0$, then let $\ell\in \omega$ be maximal with $x_{\ell-1} =0$. We obtain $v((x_m))=\vmin(x_\ell)>0$ and thus $(x_m) \in\OO_v$. If $x_0\neq 0$, then $x_i\neq 0$ for every $i\in \omega$. Hence, $v((x_m))=0$, also showing $(x_m) \in \OO_v$.
		\end{proof}
		\begin{claim}
					\label{cl:valrings2}
			For any $n\in \omega$ we have $v_0(\OO_n^\times) = \coprod_{\gamma = -2n+1}^0 A_{-\gamma}$ and therefore
			\begin{align}
				\label{eq:vnK}
				v_nK=\left(\coprod_{\gamma>-\infty}^{0}\!A_{-\gamma}\right) / \left(\coprod_{\gamma = -2n+1}^0 A_{-\gamma}\right) = \coprod_{\gamma>-\infty}^{-2n} A_{-\gamma}.
			\end{align}
		\end{claim}
		\begin{proof}[Proof of Claim~\ref{cl:valrings2}]
		\renewcommand{\qedsymbol}{$\Diamond$}
	     Let $(x_m)\in \OO_1^\times$. Then $x_1\in K_1\setminus\{0\}$. Hence, by the inclusion of ordered fields above, we have
		$x_1\in \R\pow{\coprod_{\gamma =-1}^0 A_{-\gamma}}$ and thus $v_0((x_m))=\vmin(x_1)\in\coprod_{\gamma =-1}^0 A_{-\gamma}$. 
		By iteration of this argument, we obtain \eqref{eq:vnK} for any $n\in \omega$.
		\end{proof}
		 
		To complete our construction, set $u_{2n}=v_n$ for any $n\in \omega$. Claim~\ref{cl:vnlesshens} and Claim~\ref{cl:valrings2} yield that $u_{2n}$ is a $(2n+1)_{\leq}$-henselian (and thus $(2n)_{\leq}$-henselian) valuation with $u_{2n}K=\coprod_{\gamma>-\infty}^{-2n} A_{-\gamma}$. Now let $u_{2n+1}$ be the coarsening of $u_{2n}$ induced by the convex subgroup $A_{2n}$ of $u_{2n}K$. Since $u_{2n+1}$ is the coarsening of a $(2n+1)_{\leq}$-henselian valuation, it is itself also $(2n+1)_{\leq}$-henselian. By definition of $u_{2n+1}$, its value group is given by $u_{2n+1}K=\coprod_{\gamma>-\infty}^{-2n-1} A_{-\gamma}$, as required.
		\qed
	\end{construction}
	
	\begin{remark}
		In \Autoref{constr:main} we could have started with any ordered field $K_0$ (instead of $\R$) to obtain the required $t$-henselian non-henselian ordered field. In fact, one can also start with any field $K_0$ of characteristic $0$ (so \Autoref{cor:biljakovic} is applicable) to obtain a $t$-henselian non-henselian field (without an ordering) with the specified value groups. \qed
		
		    
	\end{remark}
	
	\begin{theorem}\thlabel{thm:secondmain}
		Let $A\neq\{0\}$ be an ordered abelian group. Then there exists an ordered henselian valued field $(L,<,v)$ such that $vL=\coprod_{-\omega}A$ and $(L,<)$ admits no non-trivial $0$-$\Lor$-definable henselian valuation. In particular, $v$ is not $0$-$\Lor$-definable.
	\end{theorem}
	
	\begin{proof}
		For any $\Log$-sentence $\sigma$, let $\sigma^*$ be the $\Lvf$-sentence such that for any valued field $(F,w)$ we have $wF\models \sigma$ if and only if $(F,w)\models \sigma^*$.\footnote{Syntactically, the sentence $\sigma^*$ is obtained from $\sigma$ as follows: Any instance of quantification $Qx$ (where $Q\in \{\exists,\forall\}$ and $x$ is a variable) is replaced by $Q(x\neq 0)$, any occurence of the constant symbol $0$ is replaced by $v(1)$, any occurence of a variable $x$ is replaced by $v(x)$ and finally any expression of the form $v(t)+v(s)$ is replaced by $v(ts)$.}
		Let $G=\coprod_{-\omega}A$ and let $\Th(G)$ be the complete $\Log$-theory of $G$. We set $\Sigma$ to be the $\Lvf$-theory $\{\sigma^*\mid \sigma \in \Th(G)\}$.
		
		Let $(K,<)$ and $ u_n$ (for any $n\in \N$) be as in \Autoref{constr:main} {  for the case $A_\gamma=A$ for any $\gamma\in \omega$}.  
		For any $n\in \N$ we let $h_n$ be the $\Lvf$-sentence stating that $v$ is a non-trivial  $n_\leq$-henselian valuation. Finally, let $\Th(K,<)$ be the complete $\Lor$-theory of $(K,<)$.
		
		For any $n\in \N$, the ordered valued field $(K,<, u_n)$ satisfies $\{h_i\mid i\leq n\}$. Since $ u_nK=G$, we also have $(K, u_n)\models \Sigma$. Hence, any finite subset of the $\Lovf$-theory $\Sigma \cup \{h_n\mid n\in \N\}\cup\Th(K,<)$ is satisfiable. By the Compactness Theorem, we obtain an ordered valued field $(F,<,w)$ such that $(F,<)\models \Th(K,<)$, $wF\models \Th(G)$ and $w$ is non-trivial, $n_\leq$-henselian for any $n\in \N$ and thus henselian. We set $L=Fw\pow{G}$ and $v=\vmin$ on $L$. Since $Fw=Lv$ and $wF\equiv G=vL$, the Ax--Kochen--Ershov principle for ordered fields (cf.~\cite[Corollary~4.2]{farre}) implies	$(L,<,v)\equiv (F,<,w)$. As $(F,<)\equiv(K,<)$, we also obtain $(L,<)\equiv (K,<)$.
		
		By elementary equivalence, any non-trivial $0$-$\Lor$-definable henselian valuation in $(L,<)$ corresponds to a non-trivial $0$-$\Lor$-definable henselian valuation in $(K,\allowbreak <)$. However, since $(K,<)$ does not admit any non-trivial henselian valuation, this shows the required conclusion.
	\end{proof}
	
	\begin{remark}
		\begin{enumerate}
			\item  {Let $A\neq\{0\}$ be an ordered abelian group. In order to obtain an ordered henselian valued field $(L,<,v)$ such that $vL = \coprod_{-\omega} A$ and $v$ is not $0$-$\Lor$-definable, one can also proceed as follows:\footnote{We thank the referee for this argument and for the permission to include it in our paper.} Let $L = K\pow{\coprod_{\Z} A}$ for some ordered field $(K,<)$. Now let $v$ be the coarsening of $\vmin$ induced by the convex subgroup $\coprod_{\N} A$ of $\vmin L$. Then $Lv = K\pow{\coprod_\N A}$. Now let $w$ be the strict coarsening of $v$ induced by the convex subgroup $\coprod_{\{0\}} A$ of $vL$. Then $wL = \coprod_{-\N} A \cong vL$ and $Lw = K\pow{\coprod_\omega A} \cong Lv$. We fix $a \in \OO_w \setminus \OO_v$. Assume, for a contradiction, that there is a parameter-free $\Lor$-formula $\varphi(x)$ defining $\OO_v$. Then we have that $$(L,<) \models \neg \varphi(a).$$ By the Ax--Kochen--Ershov principle for ordered fields we have that $(L,\allowbreak<,v) \equiv (L,<,w)$, whence $\varphi(x)$ also defines $\OO_w$. Since $a\in \OO_w$, we obtain $$(L,<) \models \varphi(a),$$ a contradiction. }
			
			\item Let $(L,<,v)$ be as in \Autoref{thm:secondmain}. Since any non-trivial coarsening $w$ of $v$ is also henselian, also $w$ is not $0$-$\Lor$-definable. Moreover, for any convex subgroup $B\subseteq A$, the corresponding coarsening $w$ of $v$ with value group $\left(\coprod_{-\N}A\right)\amalg (A/B)$ is not $0$-$\Lor$-definable.

			\item 
			\Autoref{thm:secondmain} stands in contrast to \cite[Theorem~4]{hong}, which implies that any henselian valuation $v$ such that $vK$ is non-divisible and archimedean is already $0$-$\Lr$-definable. In this regard, note that the rank $-\omega$ of $\coprod_{-\omega}A$ has no minimum, whence the valuation $v$ in \Autoref{thm:secondmain} has no coarsening with archimedean value group.\qed
		\end{enumerate}
	\end{remark}
	
	We can now use \Autoref{thm:secondmain} to show that the conclusion in \Autoref{fact:hongdiscrete} cannot be strengthened to parameter-free definability.
	
	\begin{example}\thlabel{ex:counter1}
		The ordered abelian group $\coprod_{-\omega}\Z$ is discretely ordered (with $\one_0$ as least positive element). By \Autoref{thm:secondmain}, there exists an  ordered henselian valued field $(L,<,v)$ with $vL=\coprod_{-\omega}\Z$ such that $v$ is not $0$-$\Lor$-definable.\qed
	\end{example}
	
	Similarly, we show that \Autoref{fact:hongprendible}, \Autoref{fact:kkl} and \Autoref{thm:strenhong} cannot be strength\-ened to parameter-free definability.
	
	\begin{example}\thlabel{ex:counter2}
		Consider the densely ordered abelian group $$A=\setbr{\left.\!\frac{a}{2^n}\ \right|\ a,n\in \Z}\subseteq \Q.$$
		Then $A$ is regular but not $3$-divisible, and is (isomorphic to) a convex subgroup of $G:=\coprod_{-\omega}A$. \Autoref{prop:prendibleclosed} yields that $G$ is not closed in $\div{G}$. By \Autoref{thm:secondmain}, there exists an  ordered field $(L,<)$ and a henselian valuation $v$ on $L$ such that $vL=G$ and $v$ is not $0$-$\Lor$-definable.\qed
	\end{example}
	
	\begin{remark} 
	\label{rem:referee}
	    By \cite[Lemma~2.3.7]{hong2}, the following holds:\footnote{A similar argument can also be found in \cite[Proof of Proposition~3.7]{jahnke5}.} 
	           \emph{Let $(K,v)$ be a henselian valued field and let $p\in\N$ be prime. Suppose that $vK$ is densely ordered and $p$-regular but not $p$-divisible. Then $v$ is $0$-$\Lr$-definable.}

	   In contrast to this result, \Autoref{ex:counter2} shows that the mere {existence} of a convex $p$-regular subgroup that is not $p$-divisible does, in general, not ensure definability without parameters.
	\end{remark}
	
	\section{Further work}\label{sec:questions}
	
	We conclude this work by collecting open questions which may serve as a possible starting point for future work on topics related to this paper.
	
	For a class $\mathcal{C}$ of non-trivial ordered abelian groups we consider the following:
	
	\begin{itemize}
	    \item[($\dagger$)$_0$] If $(K,v)$ is a henselian valued field such that $vK \in \mathcal{C}$, then $v$ is $0$-$\Lr$-definable.
	    \\[-0.3cm]
	    \item[($\dagger$)\phantom{$_0$}] If $(K,v)$ is a henselian valued field such that $vK \in \mathcal{C}$, then $v$ is $\Lr$-definable.
	\end{itemize}

	The corresponding properties ($\cancel{\dagger}$){$_0$} and ($\cancel{\dagger}$) stand for the negations of the above; that is, a class $\mathcal{C}$ of non-trivial ordered abelian groups satisfies ($\cancel{\dagger}$){$_0$} (respectively ($\cancel{\dagger}$)) if there exists a henselian valued field $(K,v)$ such that $vK\in \mathcal{C}$ but $v$ is not $0$-$\Lr$-definable (respectively $\Lr$-definable).
	
	To summarise the results from the literature as well as this work, we distinguish between  discretely and densely ordered abelian groups.
	
	\subsection*{Discrete case}
	We consider the following two elementary classes of ordered abelian groups.
	\begin{align*}
		\mathcal{C}^{\mathrm{discr}}_0\colon & \text{all regular discretely ordered abelian groups;}\\
		\mathcal{C}^{\mathrm{discr}}_1\colon & \text{all discretely ordered abelian groups.}
	\end{align*}
	The class $\mathcal{C}^{\mathrm{discr}}_0$ consists of all ordered abelian groups that are elementarily equivalent to $\Z$ as $\Log$-structures. These are also called $\Z$-groups (cf.~e.g.~\cite[Fact~3.4]{krapp}). By \cite[Theorem~4]{hong}, $\mathcal{C}^{\mathrm{discr}}_0$ satisfies ($\dagger$)$_0$. 
	
	The ordered abelian group $\Z\amalg \Z$ serves as an example of a non-regular discretely ordered abelian group, showing that $\mathcal{C}^{\mathrm{discr}}_0$ is strictly contained in $\mathcal{C}^{\mathrm{discr}}_1$. By \Autoref{fact:hongdiscrete} the class $\mathcal{C}^{\mathrm{discr}}_1$ satisfies ($\dagger$), and due to \Autoref{ex:counter1} this class does not satisfy ($\dagger$){$_0$}.
	We thus obtain the following picture:
	
	$$\underbrace{\mathcal{C}^{\mathrm{discr}}_0}_{(\dagger)_0}\quad  \subsetneq \underbrace{\mathcal{C}^{\mathrm{discr}}_1}_{(\dagger)\text{ and } (\cancel{\dagger})_0}$$
	
	\subsection*{Dense case}
	
		We consider the following four elementary classes of ordered abelian groups $\mathcal{C}^{\mathrm{dense}}_i$ consisting of all densely ordered abelian groups $G$ satisfying the condition $c_i$, where $c_i$ is given as follows:

    \begin{itemize}
        \item [$c_0\colon$] for some prime $p$, the group $G$ is $p$-regular but not $p$-divisible;
        \item [$c_1\colon$]  for some prime $p$, some convex subgroup of $G$ is $p$-regular but not $p$-divisible;
        \item [$c_2\colon$] $G$ is not closed in $\div{G}$;
        \item [$c_3\colon$] for some prime $p$, the group $G$ contains no non-trivial $p$-divisible convex subgroup.
    \end{itemize}	
    
	
	While $\mathcal{C}^{\mathrm{dense}}_0$ satisfies ($\dagger$)$_0$, the class $\mathcal{C}^{\mathrm{dense}}_1$ satisfies ($\dagger$) but not ($\dagger$)$_0$ (see \Autoref{ex:counter2} and \Autoref{rem:referee}). By \Autoref{prop:prendibleclosed} and \Autoref{ex:counterhong}, the class $\mathcal{C}^{\mathrm{dense}}_1$ is strictly contained in $\mathcal{C}^{\mathrm{dense}}_2$. Moreover, \Autoref{thm:strenhong} yields that $\mathcal{C}^{\mathrm{dense}}_2$ satisfies ($\dagger$). 
	
	Further, $\mathcal{C}^{\mathrm{dense}}_2$ is strictly contained in $\mathcal{C}^{\mathrm{dense}}_3$ (cf.~\cite[Proposition~3.18 and Example~3.20]{krapp}). The remarks in \cite[page 1147]{delon} show that for any densely ordered abelian group $G$ not satisfying condition $c_3$ above there exists an (almost real closed) valued field $(K,v)$ with $vK=G$ such that $v$ is not $\Lr$-definable. Hence, any class of densely ordered abelian groups $\mathcal{D}$ that strictly contains $\mathcal{C}^{\mathrm{dense}}_3$ does not satisfy ($\dagger$). 
	
	As a result of the discussion above, we obtain the following picture:
	
	$$\underbrace{\mathcal{C}^{\mathrm{dense}}_0}_{(\dagger)_0}  \subsetneq \underbrace{\mathcal{C}^{\mathrm{dense}}_1\subsetneq \mathcal{C}^{\mathrm{dense}}_2}_{(\dagger)\text{ and } (\cancel{\dagger})_0}\subsetneq \mathcal{C}^{\mathrm{dense}}_3 \subsetneq \underbrace{\mathcal{D}}_{(\cancel{\dagger})}$$

	
	

	\subsection*{Open questions}
     {While from the observations above it is clear that} \textcolor{black}{the class $\mathcal{G}_0$ of ordered abelian groups for which $(\dagger)_0$ holds is a \emph{proper} subclass of the class $\mathcal{G}$ of ordered abelian groups satisfying $(\dagger)$,}  { none of the described dividing lines are yet proven to be sharp.
    Therefore, we pose the following questions:}
	\textcolor{black}{
	\begin{question}
	\begin{enumerate}
	    \item[(i)] What is the largest class of ordered abelian groups satisfying $(\dagger)$, respectively $(\dagger)_0$?
	    \item[(ii)] Is there a characterisation of $\mathcal{G}_0$ inside $\mathcal{G}$?
	\end{enumerate}
	\end{question}}
	
	\subsection*{Concluding Remark}\label{sec:conlude}
	As noted in the introduction, in this paper we only consider conditions on the value group. There are also conditions for the definability of a given henselian valuation only depending on its residue field (cf.~e.g.~\cite[Proposition~3.1 and Corollaries~3.3~\&~3.8]{jahnke4} and \cite[Theorem~5.3~(3)]{krapp}). In particular, the analogous of $(\dagger)$ and $(\dagger)_0$ can be formulated for residue fields. Such conditions will be the subject of a further publication.

	\color{black}

\end{document}